\documentclass[11pt]{article}%

\usepackage{amssymb}
\usepackage{amsmath}
\usepackage{amsfonts}

\usepackage{amssymb,color}
\definecolor{c20}{rgb}{0.,0.7,0.}
\definecolor{c30}{rgb}{0.,0.,1.}
\definecolor{c40}{rgb}{1,0.1,0.7}
\definecolor{c50}{rgb}{1,0,0}

\newcommand{\ooe}{{\rm o}}

\newcommand{\ah}{\alpha}
\newcommand{\bt}{\beta}
\newcommand{\gm}{\gamma}

\newcommand{\rcal}{{\cal R }}

\newcommand{\eps}{\epsilon}

\newcommand{\ep}{\epsilon}

\def\E{{\mathbb E}}
\def\R{{\mathbb R}}

\def\tp{{\widetilde p}}

\newcommand{\tz}{\mbox{ as } t \to 0}
\newcommand{\ninf}{\mbox{ as } n\to\infty}

\newtheorem{proposition}{Proposition}
\newtheorem{theorem}{Theorem}

\newcommand{\beq}{ \begin{equation}}
\newcommand{\eeq}{ \end{equation}}
\newcommand{\beqr}{ \begin{eqnarray}}
\newcommand{\eeqr}{ \end{eqnarray}}
\newcommand{\beqrn}{ \begin{eqnarray*}}
\newcommand{\eeqrn}{ \end{eqnarray*}}

\newcommand{\bye}{ \end{document}}
\newcommand{\noi}{\noindent}

\usepackage{amssymb,color}
\definecolor{c20}{rgb}{0.,0.7,0.}
\definecolor{c30}{rgb}{0.,0.,1.}
\definecolor{c40}{rgb}{1,0.1,0.7}
\definecolor{c50}{rgb}{1,0,0}

\title{\bf  Approximation of a random process with variable smoothness}
\author{
Enkelejd Hashorva, Mikhail\ Lifshits,  Oleg\ Seleznjev
}

\date{ }

\begin{document}
\baselineskip=3.5 ex

\maketitle

\begin{abstract}
We consider the rate of piecewise constant approximation to a locally stationary process
$X(t),t\in [0,1]$, having a variable smoothness index $\ah(t)$. Assuming that $\ah(\cdot)$
attains its unique minimum at zero and satisfies
$$
    \ah(t)=\ah_0+b t^\gm+ \ooe(t^\gm) \quad \tz,
$$
we propose a method for construction of observation points
(composite dilated design) such that the integrated mean square error
\[
 \int_0^1  \E\{(X(t)-X_n(t))^2\} dt \sim
  \frac{K}{n^{\ah_0}(\log n)^{(\ah_0+1)/\gm}}  \qquad \ninf,
\]
where a piecewise constant approximation $X_n$ is based on $N(n)\sim n$ observations of $X$.
{Further, we prove that the suggested} approximation rate is optimal, and {then} show how {to find an optimal} constant $K$.
\end{abstract}

\noi{\bf Keywords:} locally stationary processes, multifractional Brownian motion, piecewise
constant approximation.

\section{Introduction}

Probabilistic models based on the locally stationary processes with variable smoothness
became recently an object of interest and a convenient tool for applications in various
areas (such as modelling of internet traffic or artificial landscapes, finance, geophysics,
biomedicine, etc) due to their flexibility. The most known representative {random process} of this class is
a multifractional Brownian motion (mBm) independently introduced in \cite{BJR} and \cite{PL}.
A more general class of $\alpha(t)$-locally stationary Gaussian processes with a variable smoothness index $\ah(t), t\in [0,1]$,  was elaborated in
\cite{DK08}. We refer to \cite{A} for a survey and to \cite{AB,ACL,FL,S08}
for studies of particular aspects of mBm.

Whenever we need to model such processes with a given accuracy, the  approximation (time discretization)
accuracy has to be evaluated.

More specifically, consider a random  process $X(t),\,t\in[0,1]$,  with finite second moment
and variable  quadratic mean smoothness (see precise definition (\ref{def1}) below).
The process $X$ is observed at $N=N(n)$
points and a piecewise constant approximation $X_n$
is built upon these observations. The approximation performance on the entire interval is
measured by integrated mean square error (IMSE) $\int_0^1 \E\{(X(t)-X_n(t))^2\} dt$.
We construct a sequence of sampling designs (i.e., sets of observation points) taking into
account the varying smoothness of $X$ such that on a class of processes, the IMSE decreases faster
when compared to conventional regular sampling designs (see, e.g., \cite{S00}) or to quasi-regular
designs, \cite{AS11}, used for approximation of locally stationary random processes and random
processes with an isolated  singularity point, respectively.

The approximation results obtained in this paper can be used in various
problems in signal processing, e.g., in optimization of compressing digitized
signals, (see, e.g., \cite{Cal02}), in numerical analysis of random functions (see,
e.g., \cite{BC92,CMR07,CL06}), in simulation studies with controlled accuracy for
functionals on realizations of random processes (see, e.g., \cite{AS08,E86}). It is
known that piecewise constant approximation gives an optimal rate for certain class
of continuous random processes  satisfying H\"older condition (see, e.g., \cite{BS99,S00}).
In this paper we develop a technique improving this rate for a certain class of
locally stationary processes with variable smoothness. The developed technique can
be generalized for more advanced approximation methods (e.g., Hermite splines)
and various classes of random processes and fields.

Some related approximation results for continuous and smooth random functions can
be found in \cite{HPS03,KP05,S96}. The book \cite{R00} contains a very detailed
survey of various random function approximation problems.

The paper is organized as follows.
In Section \ref{se:Procmed} we specify the problem setting. We recall a notion of
locally stationary process, introduce a class of piecewise constant approximation
processes, and define integrated mean square error (IMSE) as a measure of approximation
accuracy. Furthermore, we introduce a special method of {\it composite dilated sampling
designs} that suggests how to distribute the observation points sufficiently densely
located near the point of the lowest smoothness. The implementation of this design
depends on some functional and numerical parameters, and we set up a certain number of
mild assumptions about these parameters.
In Section \ref{se:Res}, our main results are stated. Namely, for a locally
stationary process with known smoothness, we consider the piecewise constant
interpolation related to dilated sampling designs (adjusted to smoothness parameters)
and find the asymptotic behavior of its approximation error. In the second part
of that section, the approximation for conventional regular and some quasi-regular
sampling designs are studied.
In Section \ref{se:opt}, the results and conjectures related to optimality of our
bounds are discussed.
Section \ref{se:Pro} contains the proofs of the statements from Section \ref{se:Res}.

%

\section{Variable smoothness random processes and
approximation methods. Basic notation}
\label{se:Procmed}


\subsection{Approximation problem setting}

Let $X=X(t), t\in [0,1]$, be an $\ah(\cdot)$-locally
stationary random process, i.e., $\E\{X(t)^2\}<\infty$ and
\beq \label{def1}
     \lim_{s\to 0}\frac{||X(t+s)-X(t)||^2}{|s|^{\ah(t)}}=c(t)\quad
     \mbox{ uniformly in  }  t\in [0,1],
\eeq
where $||Y||:=(\E Y^2)^{1/2}$,
$\ah(\cdot), c(\cdot) \in C([0,1])$ and
$ 2\ge  \ah(t)>0,\, c(t)>0$.

We assume that the following conditions hold for the function $\ah(\cdot)$ describing
smoothness of $X$:

(C1)  $\ah(\cdot)$ attains its global minimum $\ah_0:=\ah(0)$ at the unique
point $t_0=0$.

(C2)  there exist $b,\gamma>0$ such that
$$
    \ah(t)=\ah_0+b t^\gm+ \ooe(t^\gm) \quad \tz.
$$
The choice $t_0=0$ in $(C1)$ is of course just a matter of notation convenience. The results are essentially
the same for any location of the unique minimum of $\ah(\cdot)$.

Let $X$ be sampled at the distinct design points $T_n=(t_0(n),\ldots, t_{N}(n))$ (also referred as knots), where
$0=t_0(n)<t_1(n)<\dots <t_{N}(n)=1$, $ N=N(n)$. We suppress
the auxiliary integer argument $n$  for design points $t_j=t_j(n)$ and for number of points
$N=N(n)$ when doing so causes no confusion.

The corresponding piecewise constant approximation is defined by
\[
   X_n(t):=X(t_{j-1}),   \qquad t_{j-1}\le t <t_j,\ j=1,\ldots, N.
\]

\noi In this article, we consider the accuracy of approximation to
$X$ by $X_n$ with respect to the integrated mean square error (IMSE)
\[
    e_n^2=||X-X_n||_2^2:= \int_0^1 ||X(t)-X_n(t)||^2 dt.  
\]
We  describe now a construction of sampling designs $\{T_n\}$ providing the fastest decay of
$e_n^2$.

\subsection{Sampling design construction}

The construction idea is as follows. In order to achieve a
rate-optimal approximation of
$X$ by $X_n$, we introduce a sequence of \emph{dilated} sampling designs $\{T_n\}$.

Recall first that any probability density $f(t), t\in [0,1]$, generates a sequence of associated
conventional sampling designs,  (cf., e.g., \cite{SY}, \cite{BC92},  \cite{S00}) defined by
\beq \label{jn}
    \int_0^{t_j} f(t) dt= \frac{j}{n}, \quad  \;  j=0, \ldots, n,
\eeq
i.e., the corresponding sampling points are $(j/n)$-percentiles of the density $f(\cdot)$, say,
a \emph{sampling density}.

 Let $p(\cdot) $ be a probability density on $\R_+:= [0,\infty)$. In our problem, it turns out
 to be useful to \emph{dilate} the design density by replacing it with a \emph{dilated sampling density}
\beq \label{eq:hn}
    p_n(t):= d_n \, p(d_n t) ,\qquad  t\in [0,1],
\eeq
where $d_n\nearrow \infty$ is a  \emph{dilation coefficient}. Note, that formally $ p_n(\cdot)$ is not a
probability density, but
$$
 \int_0^1  p_n(t) dt= \int_0^{d_n}  p(u) du\to 1 \ninf.
$$
The idea of dilation is obvious: we wish to put more knots near the point of the worst
smoothness. The delation coefficient should be chosen accordingly to the smoothness behavior
at this critical point. In our case, $(C2)$ requires the choice
\[
   d_n:=(\log n)^{1/\gm}
\]
that will be maintained in the sequel. As in (\ref{jn}), we define the knots by
\beq \label{eq:knots2}
    \int_0^{t_j} p_n(t) dt= \frac{j}{n}.
\eeq
Further optimization of the approximation accuracy bound requires one more adjustment:
it turns out to be useful to choose the knots $t_j$ as in (\ref{eq:knots2}) using different
densities in a neighborhood of the critical point and outside of it.
We call \emph{composite} such  sampling design constructions operating differently on two
disjoint domains.
\medskip

Now we pass to the rigorous description of our sampling designs.
Let $p(u)$ and $\tp(u),u\in [0,\infty)$, be two probability densities.
Let the dilated densities $p_n(\cdot)$ be defined as in (\ref{eq:hn}). Similarly,
$\tp_n(t):= d_n \, \tp(d_n t)$.

For $0<\rho<1$, we define the \emph{composite dilated
$(p,\rho,\tp)$-designs} $T_n$ by choosing $t_j$ according to (\ref{eq:knots2})
for
\[
     0\le j\le  J(p,\rho,n)
     :=
     n\, \int_0^\rho p_n(t)dt  
     = 
     n\, \int_0^{\rho d_n} p(u)du  
     \le n.
\]
Notice that for these knots, we have $0\le t_j\le \rho$. Furthermore, we fill the interval
$(\rho,1]$ with analogous knots $t_i$ using density $\tp(\cdot)$,
\beq \label{eq:knots2t}
    \int_0^{t_i} \tp_n(t) dt= \frac{j}{n}\, ,
\eeq
where
\[
       J(\tp,\rho,n) < j\le J(\tp, 1,n),
\]
\[
i=j+J(p,\rho,n)-J(\tp,\rho,n).
\]
For these knots we clearly have $\rho< t_i\le 1$.
Note, it follows by definition that
\[
    J(p, \rho,n) \sim n \int_0^{\rho d_n} p(u) du   \sim n \qquad \ninf,
\]
and similarly, in the interval $[\rho, 1]$, the number of points does not exceed
\[
  n - J(\tp, \rho,n) \sim n \int_{\rho d_n}^\infty \tp(u) du   = \ooe(n) \qquad \ninf,
\]
that is the total number of sampling points satisfies
\beq \label{eq:knN}
   N(n)\sim J(p,\rho,n)\sim n   \ninf.
\eeq

In the sequel, we will use $(p,\rho,\tp)$-designs satisfying the following
additional assumptions on $p(\cdot)$, $\rho$, and $\tp(\cdot)$:

 (A1) The design density $p(\cdot)$ is bounded, non-increasing, and
\begin{equation} \label{eq:A1}
     p(u)\ge q_1 \exp\{  -q_2 u^\gm \}, \qquad u\ge 0,\ q_1>0,\,  \tfrac{b}{\ah_0} > q_2>0.
\end{equation}

(A2) We assume that $\tp$ is \emph{regularly varying } at infinity
with some index $r\le -1$. This means that for all $\lambda>0$,
\beq \label{eq:RV}
     \frac{\tp(\lambda u)}{\tp(u)}\rightarrow \lambda^r \quad  \mbox{ as } u \to +\infty.
\eeq
In this case, we write $\tp(\cdot) \in\rcal_{r}(+\infty)$.

 (A3) Finally, we assume that the  parameter $\rho$ is small enough.
Namely, applying  $q_2<{b}/{\ah_0}$ and using $(C2)$ we may choose
$\rho$ satisfying
\beq \label{eq:rho}
    q_2 \sup_{0\le t \le \rho} \ah(t)
    < \inf_{0\le t \le \rho} \frac{\ah(t)-\ah(0)}{t^\gm}
\eeq
and
\beq \label{eq:rho2}
    q_2 \rho^\gm<1.
\eeq

\noi  For example,  let $\alpha(t)=1+t^\gamma$. Then $(C1)$, $(C2)$ hold and  (A3) corresponds to $\rho<(1/q_2-1)^{1/\gamma}$,
where $0<q_2<1$.
\smallskip

Regularly varying densities satisfy (\ref{eq:A1}) for large $u$, thus we could
simplify the design construction by letting $p=\tp$. However, this simplified
choice does not provide an optimal constant $K$ in the main approximation error asymptotics  (\ref{eq:main})
below.

\section{{Main} results} \label{se:Res}

\subsection{Dilated approximation designs} \label{se:S3}

\noi In the following theorem, we give the {principal} result of the paper and    consider IMSE $e_n^2$ of approximation  to
$X$ by $X_n$ for the proposed sequence of  composite dilated sampling designs $T_n, n\ge 1$.
It follows from (A1) that  the following constant is finite,
$$
   K=K(c,\alpha, (p,\rho,\tp)) :=\frac{c_0}{\ah_0+1}  \int_0^\infty p(u)^{-\ah_{0}} e^{-b u^\gm} du < \infty,
$$
where $c_0:=c(0)$.

\begin{theorem}  \label{th:As} Let $X(t), t\in[0,1]$, be an $\ah(\cdot)$-locally
stationary random process such that assumptions  $(C1)$, $(C2)$ hold.
 Let $X_n$ be the piecewise constant approximations corresponding to composite
dilated $(p,\rho,\tp)$-designs $\{T_n\}$ satisfying $(A1)$-$(A3)$.
Then $N(n)\sim n$ and
\beq \label{eq:main}
   ||X-X_n||_2^2  \sim  \frac{ K}{n^{\ah_0}(\log n)^{(\ah_0+1)/\gm}}
   \sim \frac{ K}{N^{\ah_0}(\log N)^{(\ah_0+1)/\gm}}  \qquad \ninf.
\eeq
\end{theorem}

\noi {\bf Remark 1.}\ Among the assumptions of Theorem \ref{th:As}, the monotonicity
of $p(\cdot)$ is worth of a discussion. Of course, it agrees with the heuristics
to put more knots at places where the smoothness of the process is worse. However,
this assumption may be easily replaced by some mild regularity assumptions on $p(\cdot)$.
\smallskip

\noi {\bf Remark 2.}\
The following density $p^*(\cdot)$
$$
   p^*(u)= C e^{-b u^\gm/(\ah_{0}+1)},
   \quad C=\frac{b^{1/\gm}}{(\ah_0+1)^{1/\gm} \Gamma(1/\gm+1)} \, .
$$
minimizes the constant $K$ in Theorem \ref{th:As}
and generates the asymptotically optimal sequence of designs $T_n^*$. For the optimal  $T_n^*$,
$$
   K^*:=\frac{c_0}{\ah_0+1} \left(\int_0^\infty  e^{-b u^\gm/(\ah_{0}+1)} du\right)^{\ah_0+1}=
   \frac{c_0}{\ah_0+1} \,  \left(\frac{(\ah_0+1)^{1/\gm} \Gamma(1/\gm+1)}{b^{1/\gm}}\right)^{\ah_{0}+1},
$$
see, e.g., \cite{S00}.
We stress that $p^*(\cdot)$ satisfies  assumption (\ref{eq:A1}) but it is not regularly varying.
In other words, a simple design based on $\tp=p=p^*$ does not fit in theorem's assumptions.
\smallskip

\noi {\bf Remark 3.}\ The idea of considering composite designs might seem to be overcomplicated
from the first glance. However, in some sense it can not be avoided. The previous remark shows that
if we want to get the optimal constant $K$, we must handle the exponentially decreasing densities.
Assume that
\beq \label{eq:exp}
   p(u) \le q_1 \exp\{-q_2 u^\gm\}.
\eeq
If we would simplify the design by defining $t_j(n)$ as in (\ref{eq:knots2}) for the entire interval,
i.e., with $\rho=1$, then we would have
\[
   \int_0^{t_j} p_n(t) dt =\frac{j}{n},
\]
hence,
\begin{eqnarray*}
   \frac{1}{n} &=& \int_{t_j}^{t_{j+1}} p_n(t) dt
   =  \int_{t_j}^{t_{j+1}} d_n p(d_n t) dt
   \le    d_n q_1 \int_{t_j}^{t_{j+1}} \exp\{-q_2 (d_n t)^\gm\}  dt
\\
   &\le&    d_n q_1 (t_{j+1}-t_{j}) \exp\{-q_2 (d_n t_j )^\gm\}.
\end{eqnarray*}

Let $a\in (0,1)$ and $t_j\in [1-a, 1]$. Then for the length of  corresponding intervals,
we have
$$
    t_{j+1}-t_{j} \ge \frac{\exp\{q_2 (d_n t_j)^\gm \}}{nd_n q_1}
                 \ge \frac{\exp\{q_2  \log n(1-a)^\gm \}}{nd_n q_1}.
$$
If $q_2>1$ and $a$ is so small that $q_2 (1-a)^\gm >1 $, we readily obtain
$t_{j+1}-t_{j}>a$ for large $n$ which is impossible.  Therefore, for $q_2>1$ there are no sampling points $t_j$ in $[1-a, 1]$, i.e.,
clearly  $e_n^2\ge C>0$ for any $n$, i.e.,  IMSE does not
tend to zero at all.

The confusion described above may really appear in practice because $q_2>1$ is compatible
with assumption $q_2<\tfrac{b}{\ah_0}$ from  (\ref{eq:A1}) whenever $b>\ah_0$.
\bigskip

Theorem \ref{th:As} shows that for the design densities with regularly varying tails, we may define
all knots by (\ref{eq:knots2}) without leaving empty intervals as above.
However,  we can not achieve the optimal constant $K$ on this simpler way.

\noi {\bf Remark 4.}\
Actually, the choice of knots outside of $[0,\rho]$ is not relevant for approximation rate. One
can replace the knots from (\ref{eq:knots2t}) with a uniform grid of knots $t_i= i n^{-\mu}$
with appropriate $\mu<1$.

\subsection{Regular sampling designs}

The approximation algorithm investigated in  Theorem \ref{th:As} is based
upon the assumption that we know the point where $\ah(\cdot)$ attains its minimum,
as well as the index $\gm$ in $(C2)$. If for the same process neither critical point nor index $\gm$ are
known, a conventional \emph{regular design} can be used.

Let a random process $X(t), t\in [0,1]$, be an $\alpha(\cdot)$-locally stationary,
i.e., (\ref{def1}) hold. Consider now sampling designs
$T_n=\{t_j(n), \ j=0,1,\dots, n \}$  generated by a regular positive continuous
density $p(t), t\in [0,1]$,  (see, e.g., \cite{SY}, \cite{S00}) through (\ref{eq:regular1}), i.e.,
\beq \label{eq:regular1}
  \int_0^{t_j} p(t) dt= \frac{j}{n}\ , \qquad \qquad  0\le j \le  n.
\eeq
Let  the constant
$$
   K_1:=\frac{c_0}{\ah_0+1}\; \frac{\Gamma(1/\gm+1)}{p_0^{\ah_0} b^{1/\gm}}\; , \quad p_0:=p(0).
$$

\begin{theorem}  \label{th:Un} Let $X(t), t\in[0,1]$, be an $\ah(\cdot)$-locally
stationary random process and   $(C1)$, $(C2)$ hold. Let $X_n$ be the piecewise constant
approximations corresponding to
(regular) sampling designs $\{T_n\}$ generated by $p(\cdot)$. Then
$$
   ||X-X_n||_2^2   \sim  \frac{ K_1}{n^{\ah_0}(\log n)^{1/\gm}}   \qquad \ninf.
$$
\end{theorem}
\medskip


\noi {\bf Remark 5.}\
 If  the point where $\ah(\cdot)$ attains its minimum, is known but $\gm$ is unknown,
we may build the designs without dilating the density. Instead, one could use  quasi-regular sampling designs generated by a possibly unbounded
 density $p(t), t\in (0,1]$, at the singularity point $t_0=0$ (cf., \cite{AS11}).
 For example, if $p(\cdot)$
is a density on $(0,1]$ such that
\[
  p(t) \sim A\, t^{-\kappa} \qquad \textrm{ as } t\searrow 0, \qquad 0<\kappa<1,
\]
and $t_j(n)$ are chosen through (\ref{eq:regular1}),
then for a locally stationary process $X$ satisfying $(C1)$ and $(C2)$, it is possible to show
a slightly weaker asymptotics than that of Theorem~\ref{th:As}
\[
   e_n^2\sim  \frac{K_2}{n^{\ah_0}(\log n)^{(1+\kappa\ah_0)/\gm}}\  \qquad \ninf,
\]
with $K_2:={c_0 A^{-\ah_0} \Gamma(1/\gm+1)}/({(\ah_0+1)b^{1/\gm}})$.

Of course, all above mentioned asymptotics differ only by a degree of logarithm
while the polynomial rate is determined by the minimal regularity index $\ah_0$.

\section{Optimality}
\label{se:opt}

\subsection{Optimality of the rate for piecewise constant approximations}
We explain here that the approximation rate $l_n^{-1}=n^{-\ah_0}d_n^{-(\ah_0+1)}$
achieved in Theorem \ref{th:As} is optimal in the class of piecewise constant
approximations for every locally stationary process satisfying $(C1)$ and $(C2)$. For a sampling design $T_n$, let the mesh size  $|T_n|:= \max\{(t_j-t_{j-1}), j=1,\ldots,n\}$.

\begin{proposition} \label{th:Opt1}
Let $X_n$ be piecewise constant approximations to a locally stationary process
$X$  satisfying $(C1)$ and $(C2)$ constructed according to designs
$\{T_n\}$ such that $N_n\sim n$  and $|T_n|\to 0 \ninf$.
Then
\beq \label{eq:opt1}
\liminf_{n\to\infty}\  l_n \, e_n^2 >0.
\eeq
\end{proposition}
%

\subsection{Optimality of the rate in a class of linear methods}

We explain here that the approximation rate $l_n^{-1}$
achieved in Theorem \ref{th:As} is optimal not only in the class of piecewise constant
approximations but in a much wider class of linear methods, --
at least for some locally stationary processes satisfying $(C1)$ and $(C2)$. The corresponding setting
is based on the notion of Gaussian approximation numbers, or $\ell$-numbers, that we
recall here.

Gaussian approximation numbers of a Gaussian random vector $X$ taking values in a normed space ${\cal X}$
are defined by
\beq
     \ell_n(X;{\cal X})^2 = \inf_{{x_1,\dots, x_{n-1}\atop \xi_1,\dots, \xi_{n-1}}}
     \E \left\{\left\|X-\sum_{j=1}^{n-1}\xi_j x_j \right\|^2_{\cal X}\right\},
\eeq
where infimum is taken over all $x_j\in {\cal X}$ and all Gaussian vectors $\xi=(\xi_1,\ldots, \xi_{n-1})\in \R^{n-1}$, see
\cite{KuL, L12}.
If  ${\cal X}$ is a Hilbert space, then
\[
   \ell_n(X; {\cal X})^2 = \sum_{j=n}^\infty  \lambda_j,
\]
where $\lambda_j$ is a decreasing sequence of eigenvalues of covariance operator of $X$.

Recall that a multifractional Brownian  motion (mBm)  with  a variable smoothness index (or fractality function) $\ah(\cdot)\in (0,2)$
introduced in \cite{BJR, PL} and studied in \cite{A, AB, ACL} is a Gaussian process defined
through its white noise representation
\[
   X(t) =\int_{-\infty}^\infty  \frac{e^{itu}-1}{|u|^{(\ah(t)+1)/2}}\, dW(u),
\]
where $W(t), t\in \R$, is a conventional  Brownian motion.
Notice that mBm is a typical example of a locally stationary process whenever $\ah(\cdot)$
is a continuous function.

In particular case of the constant fractality $\ah(t)\equiv \ah$, we obtain an
ordinary  fractional Brownian motion $B^\ah$, $\ah\in (0,2)$. For
$X=B^\ah$  considered as an element of ${\cal X}=L_2[0,1]$,
the behavior of eigenvalues $\lambda_j$ is well known, cf.\ \cite{Br}. Namely,
\[
   \lambda_j\sim c_\ah j^{-\ah-1} \qquad \mbox{ as } j \to \infty
\]
with some $c_\ah>0$ continuously depending on $\ah\in (0,2)$. It follows that
\[
   \ell_n(B^\ah; L_2[0,1])^2  \sim \ah^{-1} c_\ah\, n^{-\ah} \qquad \ninf.
\]
Hence, for all $n\ge 1$,
\[
   \ell_n(B^\ah; L_2[0,1])^2  \ge C_\ah\,  n^{-\ah}.
\]
Furthermore, since $B^\ah$ is a self-similar process, we can scale this
estimate from ${\cal X}=L_2[0,1]$
to ${\cal X}=L_2[0,r]$ with arbitrary $r>0$. An easy computation shows that
\[
   \ell_n(B^\ah; L_2[0,r])^2
   =   r^{\ah+1} \ell_n(B^\ah; L_2[0,1])^2 \ge C_\ah r^{\ah+1} n^{-\ah}.
\]
Let us now consider a multifractional Brownian motion $X$
parameterized by a fractality function $\ah(\cdot)$ satisfying $(C2)$.
For example, let
\beq \label{eq:ah_choice}
     \ah(t):=\ah_0+ b\, t^\gm, \qquad 0\le t\le 1,
\eeq
with $\ah_0,b>0$ chosen so small that $\ah_0+b<2$. This choice secures the necessary condition
\[
   0<\ah(t)<2, \qquad 0\le t\le 1.
\]
Then, letting $r=r_n:=d_n^{-1}$, we have
\begin{eqnarray*}
   \ell_n(X; L_2[0,1])^2&\ge& \ell_n(X; L_2[0,r_n])^2
      \ge M \ell_n(B^{\ah(r_n)}; L_2[0,r_n])^2
   \\
      &\ge& M C_{\ah(r_n)} r_n^{\ah(r_n)+1} n^{-\ah(r_n)}
      = M C_{\ah(r_n)} d_n^{-\ah(r_n)-1} n^{-\ah(r_n)}
   \\
      &\ge& C l_n^{-1} (d_n n)^{\ah_0-\ah(r_n)}
      = C l_n^{-1} (d_n n)^{-br_n^{\gm}}
   \\
      &=& C l_n^{-1}   (d_n n)^{-b(\log n)^{-1}}
      \ge {\widetilde C}\, l_n^{-1},
\end{eqnarray*}
for some positive $M, C_{\ah(r_n)} , C, {\widetilde C}$.
All bounds here are obvious except for the second inequality comparing approximation rate
of multifractional Brownian motion with that of a conventional fractional Brownian motion.

We state this fact as a separate result.

\begin{proposition} \nonumber
Let $X(t), a\le t\le b,$ be a multifractional Brownian motion corresponding to a continuous
fractality function $\ah:[a,b]\to (0,2)$. Let $B^\beta$ be a fractional Brownian motion such that
$ \inf_{a\le t\le b} \ah(t)\le \beta <2$. Then there exists $M=M(\ah(\cdot),\beta)>0$
such that
\[
    \ell_n(X; L_2[a,b]) \ge  M \ell_n(B^\beta, L_2[a,b]), \qquad n\ge 1.
\]
\end{proposition}
The proof of this proposition requires the methods very different from those used in
this article. We relegate it to another publication.

Our conclusion is that a multifractional Brownian motion with fractality function
(\ref{eq:ah_choice}) provides an example of a locally stationary process satisfying assumptions
$(C1)$ and $(C2)$ such that no linear approximation method provides a better approximation rate
than $l_n^{-1}$.

\section{Proofs}  \label{se:Pro}

\noi {\em Proof of Theorem \ref{th:As}:}

\noi We represent the IMSE $e_n^2=||X(t)-X_n(t)||^2_{2}$ as the following sum
\beq \label{eq:SS}
    e_n^2 = \sum_{j=1}^{N} \int^{t_j}_{t_{j-1}}  ||X(t)-X_n(t)||^2 dt
    = \sum_{j=1}^{N} \int^{t_j}_{t_{j-1}} ||X(t)-X(t_{j-1})||^2 dt
    =: \sum_{j=1}^{N} e_{n,j}^2.
\eeq

\noi Next, for a large $U>0$, let
$$
  e_n^2=\sum_{j=1}^N e_{n,j}^2=S_1+S_2+S_3,
$$
where the sums $ S_1, S_2, S_3$ include the terms
$e_{n,j}^2 $ such that $[t_{j-1},t_j]$ belongs to $[0,U/d_n]$, $[U/d_n, \rho]$,
and $[\rho,1]$, respectively. Let $J_1$ and $J_2$ denote  the corresponding
boundaries for index $j$.
Denote by
$$
  l_n:=n^{\ah_0}d_n^{\ah_0+1}= n^{\ah_0}(\log n)^{(\ah_0+1)/\gm}.
$$
Recall that $l_n^{-1}$ is the approximation rate announced in the theorem.
We show that only $S_1$ is relevant to the asymptotics of $e_n^2$, namely,
that $l_n S_3=\ooe(1) \ninf$, while
\beq \label{eq:oU}
  \limsup_{n\to\infty}\  l_n S_2= \ooe(1) \qquad \mbox{ as } U\to \infty.
\eeq
Let $w_j:=t_j-t_{j-1}$, $u_j:=d_n t_j$ be the normalized knots  and  denote by
 $v_j:=u_j-u_{j-1}=d_n w_j$  the corresponding dilated interval lengths.

It follows by the definition of $\ah(\cdot)$-local stationarity (\ref{def1})
that for large $n$,
\beqr \label{eq:B1}
e_{n,j}^2&=&c(t_{j-1})\int_{t_{j-1}}^{t_{j}} (t-t_{j-1})^{\ah(t_{j-1})} dt\;(1+r_{n,j})
\nonumber \\
&=&
B_{j-1}\;(t_j-t_{j-1})^{\ah(t_{j-1})+1}\;(1+r_{n,j})
\nonumber
\\
&=&
 B_{j-1}\;( v_j/d_n)^{\ah(t_{j-1})+1}\;(1+r_{n,j}),
\eeqr
where $|T_n|=\max_j w_j=\ooe(1)$ and $\max_j r_{n,j}=\ooe(1) \ninf$ and
\beqrn
B_{j}:=\frac{c(t_{j})}{\ah(t_{j})+1},\qquad j=1, \ldots, N.
\eeqrn

First, we evaluate $S_3$. Recall that for $j>J_2$ we have
$\rho d_n\le u_{j-1}<u_j\le d_n$.
We use now the following property of regularly varying functions
(see, e.g., \cite{BGT}): convergence in (\ref{eq:RV}) is uniform for all intervals
$0<a\le \lambda \le b<\infty$.
Using this uniformity we obtain, for some $C_1>0$,
\[
    \inf_{u_{j-1}\le u\le u_j}\tp(u)
    \ge \inf_{\rho d_n \le u\le d_n}\tp(u) \ge C_1\,  \tp(d_n).
\]
It follows by (\ref{eq:knots2t}) that
\[
  \int_{u_{j-1}}^{u_j} \tp(u) du =\int_{t_{j-1}}^{t_j} \tp_n(t) dt = \frac 1n\, .
\]
Hence, for some $C_2>0$,
\begin{eqnarray}
\nonumber
   v_j &\le& \left( \inf_{u_{j-1}\le u\le u_j}\tp(u)\right)^{-1} \int_{u_{j-1}}^{u_j} \tp(u) du
   \le    \frac{1}{n}\left( \inf_{u_{j-1}\le u\le u_j} \tp(u)\right)^{-1}
\\   \label{eq:V}
   &\le&  \frac{1}{C_1\, n \tp(d_n)}
   \le  C_2 \frac{ d_n^{|r|+1}}{n},\quad  j= J_2+1,\ldots, N,
\end{eqnarray}
and $\max_{j>J_2} w_j= d_n^{|r|}/{n} $.
Recall that by assumption $(C1)$,
$$
    \ah_1:=\inf_{t\in  [\rho,1]} \ah(t)>\ah_0.
$$
 Therefore, for large $n$, we get by (\ref{eq:B1}) and (\ref{eq:V}),  $\; C_3,C_4>0$,
\beq \label{eq:uppS3}
   S_3 \le n \max_{j>J_2}e_{n,j}^2\le n C_3 (v_j/d_n)^{\ah_1+1}
   \le C_4\frac{d_n^{|r|(\ah_1+1)}}{n^{\ah_1}}= \ooe(l_n^{-1}) \ninf.
\eeq

Now consider the first two zones  corresponding to $S_1, S_2$. We have by definition
$$
   \int_0^{u_j} p(u) du= \frac{j}{n}\ ,\qquad \quad  0\le j < J.
$$
Since the function $p_n(t), t\in [0,1]$,  is non-increasing, the sequence $\{v_j\}$ is non-decreasing. In fact,
\[
   \frac 1n =   \int_{u_{j-1}}^{u_j} p(u) du  \in  [p(u_j) v_j, p(u_{j-1}) v_j],
\]
and therefore,
\beq \label{eq:W}
    \frac{1} {n  p(u_{j-1})} \le v_j \le \frac{1}{n  p(u_{j})}   \le v_{j+1}\
\eeq
and it follows by (A1) that $\max_{j\le J_2} w_j =\ooe(1) \ninf$.
For $j\le J_2$,  the bounds    (\ref{eq:B1}) and (\ref{eq:W}) yield for $n$ large,
\begin{eqnarray}
  \nonumber
  e_{n,j}^2
   &=& B_{j-1} (v_j/d_n)^{\ah(t_{j-1})}\frac{v_j}{d_n} (1+\ooe(1))
  \\ \nonumber
   &\le& B_{j-1} (nd_n p(u_j))^{-\ah(t_{j-1})} \ \frac{v_j}{d_n} (1+\ooe(1))
 \\ \nonumber
   &\le& B_{j-1} (np(u_j))^{-\ah(t_{j-1})} d_n^{-\ah_0-1}   v_j (1+\ooe(1))
 \\ \label{eq:upp1}
   &=& B_{j-1} l_n^{-1} n^{ -(\ah(t_{j-1})-\ah_0) } p(u_j)^{-\ah(t_{j-1})} v_j (1+\ooe(1)).
 \end{eqnarray}
 From now on, we proceed differently in the first and in the second zone.

For the second zone, $J_1\le j\le J_2$, we do not care about the constant by using
\beq \label{eq:bstar}
  B_j\le B_*:= \max_{0\le t\le 1}\  \frac{c(t)}{\ah(t)+1}\ .
\eeq
Next, (\ref{eq:A1}) and (\ref{eq:rho}) give
\beq \label{eq:upp4}
p(u_j)^{-\ah(t_{j-1})} \le C \exp \{ q_2 \ah(t_{j-1}) u_{j}^\gm\}
\le C \exp \{ \bt_1 u_{j}^\gm\}, \quad C>0,
\eeq
where $\bt_1:= q_2 \sup_{0\le t \le \rho} \ah(t)$.
On the other hand, we infer from (\ref{eq:rho}) that
\beq \label{eq:upp3}
n^{ -(\ah(t_{j-1})-\ah_0) }
= n^{ -  \frac{\ah(t_{j-1})-\ah_0}{t_{j-1}^\gm}\, t_{j-1}^\gm }
\le  n^{ -  \bt_2 \, \frac{u_{j-1}^\gm}{ \log n } }
=  \exp\{ -  \bt_2 \, u_{j-1}^\gm \},
\eeq
where $\beta_2:= \inf_{0\le t \le \rho}({\ah(t)-\ah_0})/{t^\gm}>\bt_1$ by
(\ref{eq:rho}).

Recall that by (\ref{eq:rho2}), we have $1-q_2\rho^\gm>0$.
Moreover, for $U\le u_j\le \rho d_n$, we derive from  (\ref{eq:A1}) and (\ref{eq:W})
\[
   v_j\le n^{-1} p(\rho d_n)^{-1}
   \le   C n^{-1} \exp\{ q_2 (\rho d_n)^\gm  \}
   =  C n^{-(1-q_2\rho^\gm)},  \quad C>0,
\]
and it follows
\beq \label{eq:upp5}
   u_{j+1}^\gm- u_{j-1}^\gm
   = u_{j-1}^\gm\left( \left( \frac{u_{j+1}}{u_{j-1}} \right)^\gm-1\right)
   =  O\left( d_n^\gm n^{-(1-q_2\rho^\gm)}\right)  =\ooe(1) \ninf
\eeq
uniformly in $J_1\le j\le J_2$.

Since
 $\{v_j\}$ is non-decreasing, (\ref{eq:upp5})  implies  an integral bound
\begin{eqnarray} \nonumber
   &&     \exp\{ -  \bt_2 \, u_{j-1}^\gm \} \exp \{ \bt_1 u_{j}^\gm\} v_j
   \\ \nonumber
   &=&
   \exp\{ \bt_2 [u_{j+1}^\gm- u_{j-1}^\gm] \}
   \exp\{ \bt_1 u_{j}^\gm- \bt_2u_{j+1}^\gm \}  v_j
   \\ \nonumber
    &\le&
   C
   \inf_{ u_{j}\le u\le u_{j+1}} \exp\{ \bt_1 u^\gm- \bt_2u^\gm \} v_{j+1}
   \\   \label{eq:integral2}
   &\le&
   C
   \int^{u_{j+1}}_{u_{j}} e^{-(\bt_2-\bt_1)u^\gm} du, \quad C>0.
\end{eqnarray}

By plugging (\ref{eq:bstar}), (\ref{eq:upp3}), and  (\ref{eq:integral2}) into (\ref{eq:upp1}),
and summing up the resulting bounds over $J_1<j\le J_2$, we obtain
\beq \label{eq:uppS2}
   S_2 \le
    B_*  l_n^{-1} \int_U^\infty  e^{-(\bt_2-\bt_1)u^\gm} du\, (1+\ooe(1))  \qquad  \ninf.
\eeq
Therefore, (\ref{eq:oU}) is valid.

 In the first zone, $j\le J_1$, $t_j\le U/d_n$, the knots are uniformly small. Hence,
 $B_{j-1}$ are uniformly close to $B$ due to the continuity of the functions $\ah(\cdot)$ and
 $c(\cdot)$. Moreover, by $(C2)$ for any $\eps>0$, we have for all $n$ large enough
 \beq \label{eq:LS}
    \ah_0+(b-\ep)t_{j-1}^\gm \le \ah(t_{j-1})\le \ah_0+(b+\ep)t_{j-1}^\gm , \qquad j\le J_1 .
 \eeq
Hence (\ref{eq:upp1}) yields
 \begin{eqnarray}
  \nonumber
   e_{n,j}^2 &\le& (B+\eps) l_n^{-1} n^{-(b-\ep)t_{j-1}^\gm}\
                   p(u_j)^{-\ah_0} p(u_j)^{-(\ah(t_{j-1})-\ah_0)} \ v_j
\\ \label{eq:upp2}
             &=&   (B+\eps) l_n^{-1}\ n^{-(b-\ep)(u_{j-1}/d_n)^\gm}\
                     p(u_j)^{-\ah_0} p(u_j)^{-(\ah(t_{j-1})-\ah_0)}  \ v_j.
\end{eqnarray}
Recall that by the definition of $d_n$, we have
\[
    n^{-(b-\ep)(u_{j-1}/d_n)^\gm}
    = n^{-(b-\ep)u_{j-1}^\gm/ \log n }=\exp\{-(b-\ep)u_{j-1}^\gm \}.
\]
Since $p(\cdot)$ in non-increasing and $\{v_j\}$ is non-decreasing, we also have
an integral bound
\begin{eqnarray} \nonumber
   &&  \exp\{-(b-\ep)u_{j-1}^\gm \} \   p(u_j)^{-\ah_0}  v_j
   \\
   \nonumber
   &=&   \exp\{(b-\ep)[u_{j+1}^\gm- u_{j-1}^\gm]\}
          \ \   p(u_j)^{-\ah_0}   \exp\{-(b-\ep)u_{j+1}^\gm\} v_j
   \\
   \nonumber
   &\le&
   \exp\{(b-\ep)[u_{j+1}^\gm- u_{j-1}^\gm] \}
   \inf_{u_{j}\le u\le u_{j+1}}(p({u})^{-\ah_0} e^{-(b-\eps)u^\gm)} \ v_j
   \\ \label{eq:integral}
   &\le&
   \exp\{(b-\ep)[u_{j+1}^\gm- u_{j-1}^\gm] \}
   \int^{u_{j+1}}_{u_{j}}p({u})^{-\ah_0} e^{-(b-\eps)u^\gm} du.
\end{eqnarray}
Moreover, for $u_j\le U$, we derive from  (A1) and (\ref{eq:W})
\[
   v_j\le n^{-1} p(U)^{-1}.
\]
By using convexity  and   concavity of the power function for $\gm\ge 1$
and  $\gm\le 1$, respectively,
we get
\begin{eqnarray*}
   u_{j+1}^\gm- u_{j-1}^\gm
   &\le& \gm U^{\gm-1} (u_{j+1}-u_j)
   = \gm U^{\gm-1} (v_j+v_{j+1})
\\
   &\le& 2\gm U^{\gm-1} v_{j+1}=\ooe(1) \qquad \ninf \qquad (\gm\ge 1);
\\
   u_{j+1}^\gm- u_{j-1}^\gm
   &\le& (u_{j+1}-u_{j-1})^\gm
\\
   &=&  (v_j+v_{j+1})^\gm =\ooe(1) \qquad \ninf \qquad (\gm\le 1).
\end{eqnarray*}
Therefore, the exponential factor in (\ref{eq:integral}) turns out to be negligible.

Finally, for  $u_j\le U$, the property $d_n\to\infty$ yields
\beq \label{eq:p1}
    p(u_j)^{-(\ah(t_{j-1})-\ah_0)}  \le   \max\{ 1, p(U)^{ -  \max_{0\le t\le U/d_n} (\ah(t)-\ah_0)}\}
    = 1+\ooe(1).
\eeq
By plugging (\ref{eq:integral}) and (\ref{eq:p1}) into (\ref{eq:upp2}),
and summing up the resulting bounds over $j\le J_1$, we obtain
\[
   S_1 \le
    (B+2\eps)  l_n^{-1} \int_0^\infty p({u})^{-\ah_0} e^{-(b-\eps)u^\gm} du \
       \qquad  \ninf.
\]
Since $\eps$ can be chosen arbitrarily small, we arrive at
\beq  \label{eq:uppS1}
  \limsup_{n\to\infty}\  l_n S_1 \le  B \int_0^\infty p({u})^{-\ah_0} e^{-bu^\gm} du
  = K.
\eeq

Combining  (\ref{eq:uppS3}), (\ref{eq:uppS2}), and (\ref{eq:uppS1}) gives the desired upper bound.
\medskip

The lower bound is obtained along the same lines: we neglect $S_2$ and $S_3$,
and evaluate $S_1$ starting again from (\ref{eq:B1}). As in (\ref{eq:upp1}), we have
\begin{eqnarray}
  \nonumber
  e_{n,j}^2
   &=& B_{j-1} (v_j/d_n)^{\ah(t_{j-1})}\frac{v_j}{d_n} (1+\ooe(1))
  \\ \nonumber
   &\ge& B_{j-1} (nd_np(u_{j-1}))^{-\ah(t_{j-1})} \ \frac{v_j}{d_n} (1+\ooe(1))
 \\ \nonumber
   &=& B_{j-1} n^{-\ah_0} n^{ -(\ah(t_{j-1})-\ah_0) } p(u_{j-1})^{-\ah(t_{j-1})} d_n^{-\ah_0-1} d_n^{\ah_0-\ah(t_{j-1})}
   v_j (1+\ooe(1))
 \\ \label{eq:low1}
   &=& B_{j-1} l_n^{-1} n^{ -(\ah(t_{j-1})-\ah_0) } p(u_{j-1})^{-\ah(t_{j-1})} d_n^{\ah_0-\ah(t_{j-1})}  v_j (1+\ooe(1)).
 \end{eqnarray}
 Recall that for $j\le J_1$, coefficients
 $B_{j-1}$ are uniformly close to $B$. Moreover, by using (\ref{eq:LS}), we have for large $n$,
 \[
    d_n^{\ah_0-\ah(t_{j-1})} \ge  d_n^{- (b+\ep)t_{j-1}^\gm}
    \ge  d_n^{- (b+\ep)(U/d_n)^\gm}= 1+\ooe(1).
 \]
 Hence, (\ref{eq:low1}) yields
 \begin{eqnarray}
  \nonumber
   e_{n,j}^2 &\ge& (B-\eps) l_n^{-1} n^{-(b+\ep) t_{j-1}^\gm}\
                      p(u_{j-1})^{-\ah_0} p(u_{j-1})^{-(\ah(t_{j-1})-\ah_0)}   \ v_j
\\ \label{eq:low2}
             &=&   (B-\eps) l_n^{-1}\ n^{-(b+\ep)(u_{j-1}/d_n)^\gm}\
                     p(u_{j-1})^{-\ah_0} p(u_{j-1})^{-(\ah(t_{j-1})-\ah_0)}    \ v_j,
\end{eqnarray}
where as before
\[
    n^{-(b+\ep)(u_{j-1}/d_n)^\gm}
    = n^{-(b+\ep)u_{j-1}^\gm/ \log n }=\exp\{-(b+\ep)u_{j-1}^\gm \}.
\]
Since $p(\cdot)$ is non-increasing and $\{v_j\}$ is non-decreasing, we also have
an integral bound
\begin{eqnarray} \nonumber
   &&  \exp\{-(b+\ep)u_{j-1}^\gm \}  p(u_{j-1})^{-\ah_0}  v_j
   \\
   \nonumber
   &=&   \exp\{(b+\ep) [u_{j-2}^\gm- u_{j-1}^\gm]\}
             p(u_{j-1})^{-\ah_0}   \exp\{-(b+\ep)u_{j-2}^\gm\} v_j
   \\
   \nonumber
   &\ge&
   \exp\{(b+\ep)[u_{j-2}^\gm- u_{j-1}^\gm] \}
   \inf_{u_{j-2}\le u\le u_{j-1}} (p({u})^{-\ah_0} e^{-(b+\eps)u^\gm}) \ v_{j-1}
   \\ \label{eq:integral3}
   &\ge&
   \exp\{(b+\ep)[u_{j+1}^\gm- u_{j-1}^\gm] \}
   \int_{u_{j-2}}^{u_{j-1}}  p({u})^{-\ah_0} e^{-(b+\eps)u^\gm} du.
\end{eqnarray}
We have already seen that the exponential factor in (\ref{eq:integral3}) is negligible.

Finally, for  $u_j\le U$, the property $d_n\to\infty$ implies (cf., (\ref{eq:p1})
\beq \label{eq:p1l}
   p(u_{j-1})^{-(\ah(t_{j-1})-\ah_0)}
   \ge   \min\{ 1, p(0)^{ -  \max_{0\le t\le U/d_n} (\ah(t)-\ah_0)}\}
   = 1+\ooe(1).
\eeq
By plugging (\ref{eq:integral3}) and  (\ref{eq:p1l}) into (\ref{eq:low2}),
and summing up the resulting bounds over $j\le J_1$, we obtain
\[
   S_1 \ge
    (B-2\eps)  l_n^{-1} \int_0^U p({u})^{-\ah_0} e^{-(b+\eps)u^\gm} du \
     \qquad  \ninf.
\]
Since $\eps$ can be chosen arbitrarily small, we arrive to
\[
    \liminf_{n\to\infty}\  l_n S_1 \ge  B \int_0^U p({u})^{-\ah_0} e^{-bu^\gm} du.
\]
Finally,
\beq  \label{eq:lowS1}
 \liminf_{n\to\infty}\  l_n e_n^2
 \ge   \sup_{U>0}\  \liminf_{n\to\infty}\  l_n S_1
 \ge  B \int_0^\infty p({u})^{-\ah_0} e^{-bu^\gm} du =K.
\eeq
This is the desired lower bound.
%
\hfill $\Box$
\bigskip

\noi {\em Proof of Theorem \ref{th:Un}:}

 \noi Applying the notation of Theorem \ref{th:As}, we have for an interval approximation error
\[
    e_{n,j}^2
    = B_{j-1}\;w_j^{\ah(t_{j-1})+1}\;(1+r_{n,j}),\quad w_j= t_j-t_{j-1}, j=1,\ldots, n,
\]
where $\max_j r_{n,j}=\ooe(1) \ninf$.
Now  for a small enough $\rho>0$, similarly to Theorem~\ref{th:As}, we get
\beqrn
     \int_{\rho}^1e_n(t)^2 dt\le C/n^{\ah_1},\; C>0,\qquad
     \ah_1:=\inf_{t\in  [\rho,1]} \ah(t)>\ah_0,
\eeqrn
that is only $e_{n,j}$ such that $[t_{j-1}, t_j]\subset [0,\rho]$ are relevant
for the asymptotics, say,  $e_{n,j}, j=1, \ldots, J=J(\rho,n)$.
Let us denote the approximation rate $ L_n:= n^{\ah_0}(\log n)^{1/\gm}$.
Next, for $ S_1:=\sum_1^J e_{n,j}$ and small enough $\rho$, we have by continuity of the density
$p(\cdot)$ and the mean value theorem
\begin{eqnarray*}
   && e_{n,j}^2= B_{j-1} \;(n p(\eta_j))^{-\ah(t_{j-1})} w_j \;(1+\ooe(1)) \nonumber
\\
   &\le& \frac{B} {p(0)^{\ah_0}} (1+\eps)\, n^{-\ah_0}
   \int_{t_{j-2}}^{t_{j-1}}  e^{-(b-\eps)t^\gm \log n} dt\; (1+\ooe(1))\nonumber
\\
   &=& L_n^{-1} \; \frac{B} {p_0^{\ah_0}} (1+\eps)
   \int_{u_{j-2}}^{u_{j-1}}  e^{-(b-\eps)u^\gm} du\; (1+\ooe(1)),
\end{eqnarray*}
where $p_0:=p(0)$.
Now by summing up, we obtain
\[
    \limsup_{n\to \infty} L_n S_1 \le   \frac{B} {p_0^{\ah_0}} (1+\eps)
     \int_{0}^{\infty}  e^{-(b-\eps)u^\gm} du
     = \frac{B} {p_0^{\ah_0}} (1+\eps) \frac{\Gamma(1/\gm+1)}{(b-\eps)^{1/\gm}}
\]
Hence,
\[
    \limsup_{n\to \infty} L_n e_n^2
    =  \limsup_{n\to \infty} L_n S_1
    \le  \frac{B} {p_0^{\ah_0}} \, (1+\eps) \, \frac{\Gamma(1/\gm+1)}{(b-\eps)^{1/\gm}}.
\]
Since $\eps$ can be chosen arbitrary small, we get
$$
      \limsup_{n\to \infty} L_n e_n^2
      \le   \frac{B} {p_0^{\ah_0}}\, \frac{\Gamma(1/\gm+1)}{b^{1/\gm}}
       = K_1 .
$$
The lower bound follows similarly. This completes the proof.  \hfill $\Box$
\bigskip

\noi {\em Proof of Proposition \ref{th:Opt1}:}

\noi   Let $r_n:=d_n^{-1}=(\log n)^{-1/\gm}$ and $J_n:=\inf\{j: t_j=t_j(n)\ge  r_n\}$.
Then  (\ref{eq:B1}) implies 
\[
   e_n^2\ge \sum_{j=1}^{J_n} e_{n,j}^2 =  \sum_{j=1}^{J_n} B_{j-1} w_j^{\ah(t_{j-1})+1} (1+\ooe(1))
   =  B  \sum_{j=1}^{J_n}  w_j^{a_n+1} (1+\ooe(1)),
\]
where $a_n:=\sup_{0\le t\le r_n} \ah(t)$ and $w_j=t_j-t_{j-1}$. By using convexity of the power
function $w\to w^{a_n+1}$, we obtain
\[
   \frac{1}{J_n}  \sum_{j=1}^{J_n}  w_j^{a_n+1}
   \ge \left( \frac{1}{J_n}  \sum_{j=1}^{J_n}  w_j\right)^{a_n+1}
   \ge \left( \frac{r_n}{J_n} \right)^{a_n+1},
\]
hence,
\[
   \sum_{j=1}^{J_n}  w_j^{a_n+1} \ge \frac{r_n^{a_n+1}}{J_n^{a_n}}
   \ge \frac{r_n^{a_n+1}}{N_n^{a_n}}\ ,
\]
whereas
 \begin{eqnarray*}
   e_n^2 &\ge&  B \ \frac{r_n^{a_n+1}} {N_n^{a_n}} \ (1+\ooe(1))
   =  B \ \frac{1}{d_n^{a_n+1}}  \ \frac{1} {N_n^{a_n}} \ (1+\ooe(1))
 \\
   &=&  B  \ \frac{1}{d_n^{\ah_0+1} n^{\ah_0}} \  \left(\frac{1}{d_n n}\right)^{a_n-\ah_0}
   \    \left(\frac{n} {N_n}\right)^{a_n} \  (1+\ooe(1))
 \\
   &=& B \ l_n^{-1}   \left(\frac{1}{d_n n}\right)^{a_n-\ah_0}  (1+\ooe(1)) .
\end{eqnarray*}
Recall that  by $(C2)$,  $ a_n-\ah_0= O(r_n^\gm)= O((\log n)^{-1})$ and thus
(\ref{eq:opt1}) follows. \hfill $\Box$
\bigskip

\section*{Acknowledgments}
Research of E. Hashorva was supported  by Swiss National Science Foundation Grant 200021-1401633/1.
Research of M. Lifshits was supported by RFBR grants 10-01-00154à and 11-01-12104-ofi-m.
Research of O. Seleznjev  was  supported by the Swedish Research Council grant 2009-4489.

\noindent{\bf Authors' Addresses:}

\vskip 1pc
\noindent
E. Hashorva, \\
Actuarial Department, Faculty HEC\\
University of Lausanne,\\
CH-1015 Lausanne, Switzerland,\\
email: Enkelejd.Hashorva@unil.ch
\smallskip

\noindent M.A. Lifshits, \\
Department of Mathematics and Mechanics,\\
St.Petersburg State University,\\
198504 St.Petersburg, Russia,\\
email: mikhail@lifshits.org
\smallskip

\noindent
O. Seleznjev,\\
Department of Mathematics and Mathematical Statistics,\\
Ume{\aa} University,\\
SE-901 87 Ume\aa, Sweden, \\
email: oleg.seleznjev@matstat.umu.se

 \end{document}